\documentclass{amsart}
\usepackage{epsf, amssymb, amsfonts}

\newtheorem*{defn}{Definition}
\newtheorem{theorem}{Theorem}
\newtheorem*{lemma}{Lemma}

\newcommand{\Z}{{\bf Z}}

\title{On the notion of lower central series for loops}

\author[J.Mostovoy]{Jacob Mostovoy}
\address{Instituto de Matem\'{a}ticas (Unidad Cuernavaca)\\
Universidad Nacional Aut\'{o}noma de M\'{e}xico\\
A.P. 273 Admon. de Correos \# 3\\
C.P. 62251, Cuernavaca, Morelos, MEXICO}
\email{jacob@matcuer.unam.mx}

\begin{document}

\maketitle

The commutator calculus is one of the basic tools in group theory.
However, its extension to the non-associative context, based on the usual
definition of the lower central series of a loop, is not entirely satisfactory.
Namely, the graded abelian group associated to the lower central series of
a loop is not known to carry any interesting algebraic structure.
In this note we construct a new generalization of the lower central series to
arbitrary loops that is tailored to produce a set of multilinear operations on 
the associated graded group. 

\subsection*{The lower central series}

Let $N$ be a normal subloop of a loop $L$. There exists a unique smallest
normal subloop $M$ of $L$ such that $N/M$ is contained in the centre of the
loop $L/M$, that is, all elements of $N/M$ commute and and associate with
all elements of $L/M$. This subloop $M$ is denoted by $[N,L]$. The lower
central series of $L$ is defined by setting $\gamma_1 L=L$ and
$\gamma_{i+1} L=[\gamma_i L, L]$ for $i\geq 1$. This definition can be found
in \cite{Bruck}. For $L$ associative it coincides with the usual definition of
the lower central series.

The terms of the lower central series of $L$ are fully invariant normal
subloops of $L$ and the successive quotients $\gamma_i L/ \gamma_{i+1} L$ are
abelian groups. If $L$ is a group, the commutator operation on $L$ induces
a bilinear operation (Lie bracket) on the associated graded group
$\oplus \gamma_i L/ \gamma_{i+1} L$; this Lie bracket is compatible with the
grading. (See Chapter 5 of \cite{MKS} for a classical account of the 
commutator calculus in groups.) In general, however, there is no 
obvious algebraic structure on $\oplus \gamma_i L/ \gamma_{i+1} L$.

\subsection*{Commutators, associators and associator deviations.}
Here we introduce some terminology. 
The definitions of this paragraph are valid for $L$ a left quasigroup, that 
is, a halfquasigroup with left division.

Define the commutator of two elements $a,b$ in $L$ as
\[[a,b]:=(ba)\backslash(ab),\]
and the associator of three elements $a,b,c$ in $L$ as
\[(a,b,c):=(a(bc))\backslash((ab)c).\]
The failure of the associator to be distributive in each variable is measured
by three operations that we call {\em associator deviations} or simply {\em deviations}.
These are defined as follows:
\[(a,b,c,d)_1:=((a,c,d)(b,c,d))\backslash(ab,c,d),\]
\[(a,b,c,d)_2:=((a,b,d)(a,c,d))\backslash(a,bc,d),\]
\[(a,b,c,d)_3:=((a,b,c)(a,b,d))\backslash(a,b,cd).\]
The deviations themselves are not necessarily distributive and their failure
to be distributive is measured by the {\em deviations of the second level}. 
The general definition of a deviation of $n$th level is as follows. Given a
positive integer $n$ and an ordered set $\alpha_1, \ldots, \alpha_n$ of 
not necessarily distinct integers satisfying $1\leq\alpha_k\leq k+2$,
the deviation $(a_1,\ldots,a_{n+3})_{\alpha_1, \ldots, \alpha_n}$
is a function $L^{n+3}\to L$ defined inductively by
\[(a_1,\ldots,a_{n+3})_{\alpha_1, \ldots, \alpha_n}:=
(A(a_{\alpha_n})A(a_{\alpha_{n}+1}))\backslash
A(a_{\alpha_n}a_{\alpha_{n}+1}),\]
where $A(x)$ stands for
$(a_1,\ldots,a_{\alpha_n-1},x,a_{\alpha_n+2},
\ldots, a_{n+3})_{\alpha_1, \ldots, \alpha_{n-1}}$.
The integer $n$ is called the {\em level} of the deviation. There are 
$(n+2)!/2$ deviations of level $n$.  The associators are the deviations of level zero 
and the associator deviations are the deviations of level one.

\subsection*{The commutator-associator filtration}
Let us explain informally our approach to generalizing the lower central 
series to loops. We want to construct a filtration by normal subloops $L_i$
($i\geq 1$) on an
arbitrary loop $L$ with the following properties:
\begin{itemize}
\item[(a)] the subloops $L_i$ are fully invariant, that is, preserved by all
automorphisms of $L$;
\item[(b)] the quotients $L_i/L_{i+1}$ are abelian for all $i$;
\item[(c)] the commutator and the associator in $L$ induce well-defined
operations on the associated graded group $\oplus L_i/L_{i+1}$; these
operations should be linear in each argument and should respect the grading.
\end{itemize}
Clearly, we also want the filtration $L_i$ to coincide with the lower central
series for groups.

A na\"\i ve method of constructing such a filtration would be  setting $L_1=L$
and taking $L_i$ to be the subloop normally generated by
\begin{itemize}
\item[(a)]{commutators of the form
$[a,b]$ with $a \in L_p$ and $b \in L_{q}$ with $p+q\geq i$,}
\item[(b)] {associators of the form
$(a,b,c)$ with $a\in L_p$, $b \in L_q$, and $c\in L_r$ with
$p+q+r\geq i$.}
\end{itemize}
The subloops $L_i$ constructed in this way are fully invariant in $L$ and the
quotients $L_{i}/L_{i+1}$ are abelian groups. It can be seen that 
the commutator on $L$ induces a
bilinear operation on the associated graded group. However, the associator
does not descend to a trilinear operation on $\oplus L_{i}/L_{i+1}$.
The situation can be mended by adding to the generators of  $L_i$ for all $i$
every element of the form $(a,b,c,d)_\alpha$ with $a\in L_p$, $b\in L_q$,
$c\in L_r$ and $d\in L_s$, where $p+q+r+s\geq i$ and $1\leq \alpha \leq 3$.
This forces the associator to be trilinear on $\oplus L_{i}/L_{i+1}$, but
now we are faced with a new problem. Having introduced the deviations into the
game, we would like them to behave in some sense like commutators and
associators, namely, to induce multilinear operations that respect the grading
on $\oplus L_{i}/L_{i+1}$. This requires adding deviations of the second level
to the generators of the $L_i$, etc.

The above reasoning is summarised in the following
\begin{defn}
For a positive integer $i$, the $i$th commutator-associator subloop $L_i$ of a 
loop $L$ is $L$ itself if $i=1$, and is the subloop normally generated by
\begin{itemize}
\item[(a)] $[L_p,L_q]$ with $p+q\geqslant i$,
\item[(b)] $(L_p,L_q,L_r)$ with $p+q+r\geqslant i$,
\item[(c)] $(L_{p_1},\ldots,L_{p_n})_{\alpha_1,\ldots,\alpha_n}$
with $p_1+\ldots+p_n\geqslant i$ for all possibles choices of
$\alpha_1,\ldots,\alpha_n$.
\end{itemize}
\end{defn}

We refer to the filtration by commutator-associator subloops as the
{\em commutator-associator filtration}\footnote{for the lack of better name.}. 
In our terminology the usual
commutator-associator subloop becomes the ``second commutator-associator
subloop''. By virtue of its construction, the commutator-associator filtration 
has the desired properties: the $L_i$ are fully invariant and normal, the
quotients $L_i/L_{i+1}$ are abelian, and the associator and
the deviations induce multilinear operations on the associated graded group.
The bilinearity of the commutator is also readily seen. Indeed, take
$a,b$ in $L_p$ and $c$ in $L_q$. Modulo $L_{p+q+1}$, the commutators 
$[a,c]$ and $[b,c]$ commute and associate with all 
elements of $L$. Also, any associator that involves one element of $L_{p}$, 
one element of $L_{q}$ and any element of $L$, is trivial modulo $L_{p+q+1}$. 
Hence,
$$a(cb)\cdot[b,c]\equiv a(bc)\equiv (ab)c$$
and 
$$c(ab)\cdot [a,c]\equiv (ca)b\cdot [a,c]\equiv 
(ca\cdot[a,c])b\equiv (ac)b\equiv a(cb)$$
modulo $L_{p+q+1}$. Therefore we have
$$(c(ab)\cdot [a,c])[b,c]\equiv (ab)c \mod{L_{p+q+1}}$$
and it follows that 
$$[a,c][b,c]\equiv[ab,c]\mod{L_{p+q+1}}.$$
The linearity of the commutator in the second argument is proved in the same 
manner.

Identifying the algebraic structure on $\oplus L_i/L_{i+1}$ induced by the 
commutator, the associator and the deviations is a task beyond
the scope of the present note. It is rather easy to see that the Akivis
identity 
\begin{multline*}
[[a,b],c]+[[b,c],a]+[[c,a],b]\\
=(a,b,c)+(b,c,a)+(c,a,b)-(a,c,b)-(c,b,a)-(b,a,c)
\end{multline*}
is satisfied in $\oplus L_i/L_{i+1}$ for any $L$. One should expect that
the structure on  $\oplus L_i/L_{i+1}$ generalizes Lie rings in the same 
way as Sabinin algebras\footnote{formerly known as hyperalgebras, 
see \cite{MixSabinin} and \cite{ShestUmir}.} generalize Lie algebras.
No axiomatic definition of such a structure is yet known.

\subsection*{The lower central series and the commutator-associator subloops}
Let us compare the commutator-associator filtration with the usual lower 
central series.

It is clear that for any $L$
\[\gamma_2 L= L_2.\]
More generally, an easy inductive argument establishes that
$\gamma_i L$ is contained in $L_i$ for all $i$. There is, however, no 
converse to this statement.
\begin{theorem}\label{t1}
Let $F$ be a free loop. For any positive integer $k$ the term $F_k$ of the
commutator-associator filtration is not contained in $\gamma_3 F$.
\end{theorem}

\begin{proof}
Let $X$ be a set freely generating the loop $F$ and assume $y$ is in $X$.
Let $y^m$ stand for the right-normed product $(\ldots((yy)y)\ldots)y$ of
$m$ copies of $y$. We shall prove that for any $i\geqslant 0$ and any positive 
$m$ the deviation of $i$th level $(y^m,y,\ldots,y)_{1,1,\ldots,1}$ does not 
belong to $\gamma_3 F$. The machinery suitable for this purpose was developed 
by Higman in \cite{Higman}.

Let $L$ be a loop and assume there is a surjective homomorphism
$p: F\to L$ with kernel $N$. Define $B$ to be an additively written free
abelian group, with free generators $f(l_1,l_2)$ for all $l_1\neq 1$ and
$l_2\neq 1$ in $L$, and $g(x)$ for all $x$ in $X$. The product set
$L \times B$ can be given the structure of a loop setting
\[(l_1,b_1)(l_2,b_2)=(l_1 l_2, b_1+b_2+f(l_1,l_2)),\]
\[(l_1,b_1)/(l_2,b_2)=(l_1/l_2, b_1-b_2-f(l_1/l_2,l_2)),\]
\[(l_2,b_2)\backslash(l_1,b_1)=(l_2\backslash l_1, b_1-b_2-
f(l_2,l_2\backslash l_1)),\]
where $f(l,1)=f(1,l)=0$ for all $l\in L$. The pair $(1,0)$ is
the identity. Higman denotes this loop by $(L, B)$.

There is a homomorphism $\delta: F\to (L, B)$ defined on the generators by
\[\delta x:=(px,g(x)).\] Higman proved (Lemma~3 of \cite{Higman}) that the
kernel of $\delta$ is precisely $[N,F]$.

Without loss of generality we can assume that $F$ is the free loop on
one generator $y$; take $N=\gamma_2 F$. The quotient loop $L=F/N$ can be
identified with the group of integers $\Z$. We shall see that for any positive 
$m$ the element $\delta (y^m,y,\ldots,y)_{1,1,\ldots,1}$ is not zero in 
$(L,B)$ and, hence, that $(y^m,y,\ldots,y)_{1,1,\ldots,1}$ does not belong to 
$[N,F]=\gamma_3 F$.

\begin{lemma}
For any $m\geq 1$ and any $n\geq 0$
\[\delta (y^m,y,\ldots,y)_{\underbrace{
\scriptstyle 1,1,\ldots,1}_{\scriptstyle n}}=(0, f(n+m-1,1)+
\sum_{p<n+m-1, q} a_{p,q}f(p,q)),\]
where $a_{p,q}$ are integers.
\end{lemma}

The lemma is proved by induction on $n$. A straightforward
calculation gives
\[\delta y^m=(m,mg(y)+f(1,1)+f(2,1)+\ldots+f(m-1,1)).\]
It follows that
\[ \delta(y^{m-2}y^2)= (m,mg(y)+\sum_{p,q} c_{p,q} f(p,q))\]
where $p\leq m-2$, $q\leq 2$ and the $c_{p,q}$ are some coefficients whose
value is of no importance. Thus $\delta (y^m,y,y)$ is of the form
\[(0,f(m-1,1)+\sum_{p<m-1, q} a_{p,q}f(p,q))\] and the lemma is satisfied for
$n=0$.

Now, for arbitrary positive $n$ we have
\begin{multline*}
(y^m,y,\ldots,y)_{\underbrace{\scriptstyle 1,\ldots,1}_{\scriptstyle n}}\\
=\left.\left(
(y^m,y,\ldots,y)_{\underbrace{\scriptstyle 1,\ldots,1}_{\scriptstyle n-1}}
(y,y,\ldots,y)_{\underbrace{\scriptstyle 1,\ldots,1}_{\scriptstyle n-1}}
\right)
\right\backslash
(y^{m+1},y,\ldots,y)_{\underbrace{\scriptstyle 1,\ldots,1}_{\scriptstyle n-1}}
\end{multline*}
and it only remains to apply $\delta$ to both sides and use the induction
assumption.
\end{proof}

Higman \cite{Higman} has proved that the terms of the lower central series of 
a free loop intersect in the identity. We have no proof of a similar statement 
for the commutator-associator filtration; in view of Theorem~\ref{t1}, 
Higman's result is not sufficient to establish it.

Theorem~\ref{t1} reflects a fundamental difference between the lower central 
series  and the commutator-associator filtration. In the proof of 
Theorem~\ref{t1} we have, in fact, established that the quotient 
$\gamma_2 F/\gamma_3 F$ is not finitely generated. On the other hand, we have 
the following

\begin{theorem}
For any finitely generated loop $L$ the abelian groups $L_i/L_{i+1}$ are 
finitely generated 
for all $i>0$.
\end{theorem}
\begin{proof}
For $i=1$ the theorem is clearly true: $L_1/L_2$ is generated by the classes 
of the generators of $L$. Assume that for all $k<n$ the group $L_k/L_{k+1}$
has a finite set of generators $a_{k,\alpha}$ with $\alpha$ in some finite 
index set. Consider a commutator 
$[x,y]$ with $x$ in $L_{p}$ and $y$ in $L_{n-p}$. The class of $[x,y]$ in
$L_n/L_{n+1}$ only depends on the classes of $x$ and $y$ in $L_{p}/L_{p+1}$ 
and $L_{p}/L_{n-p+1}$ respectively. As the commutator is bilinear on 
$\oplus L_i/L_{i+1}$, the class of $[x,y]$ is a linear combination of elements
of the form $[a_{p,\alpha}, a_{n-p,\beta}]$. Similarly, all classes of 
associators and deviations can be expressed as linear combinations of 
associators and deviations of a finite number of elements. Hence,
$L_n/L_{n+1}$ is finitely generated.
\end{proof}

\subsection*{Further comments} In general, the lower central series of
loops is difficult to treat. There is one notable exception, namely,
the theory of associators in commutative Moufang loops which is analogous to
a large extent to the theory of commutators in groups, see \cite{Bruck}.
There is also another filtration, called distributor series, on
commutative Moufang loops. In general, neither of these two filtrations
coincides with the commutator-associator filtration. We note that the
deviations first appeared in the nilpotency theory of commutative Moufang
loops as the operations that produce the distributor series.

There exists a filtration on any group that is very closely related to
the lower central series, namely, the filtration by dimension subgroups.
In particular, the graded abelian groups associated to both filtrations
become isomorphic after tensoring with the field of rational numbers.
A detailed comparison of the lower central series and the dimension subgroups
can be found in \cite{Hartley}. The dimension filtration can be generalized to 
loops; it would be interesting to compare it with the commutator-associator 
filtration.

The constructions of this note can be carried out in greater generality.
The definitions of the commutator, the associator and the deviations of all
levels only involve left division and, therefore, make sense in any left
halfquasigroup. In particular, the commutator-associator filtration can be 
defined for left loops. 

\subsection*{Acknowledgment} I would like to thank Liudmila Sabinina for
many discussions.

\end{document}